\documentclass[12pt]{amsart}
\usepackage{wrapfig}
\usepackage{geometry, graphicx, latexsym, amssymb, amsmath,lscape,subfigure, bbm}
\DeclareGraphicsExtensions{.pdf}
\usepackage[all]{xy}
 \usepackage[parfill]{parskip}

\numberwithin{equation}{section}

\newcommand{\R}{\mathbb{R}}

\newcommand{\set}[1]{\left\{#1\right\}}
\newcommand{\paren}[1]{( #1 )}

\newcommand{\tabb}{\begin{center} \begin{tabular}{|l|c|c|}
\hline}
\newcommand{\tabe}{\hline \end{tabular} \end{center}}

\newcommand{\arcsec}{\text{arcsec}\,}

\begin{document}
\title{\textbf{Measuring Shape With Topology}}
\author{Robert MacPherson}
\author{Benjamin Schweinhart}
\thanks{}
\date{\today}

\begin{abstract}
We propose a measure of shape which is appropriate for the study of a complicated geometric structure, defined using the topology of neighborhoods of the structure.  One aspect of this measure gives a new notion of fractal dimension.  We demonstrate the utility and computability of this measure by applying it to branched polymers, Brownian trees, and self-avoiding random walks.
\end{abstract}

\maketitle

\section{Introduction.}\label{Introduction}

\subsection{Complicated Geometric Structures.} \label{SETS}
\begin{figure}[htp]
\centering
\subfigure[Branched Polymer]{\label{fig:BP}\includegraphics[width=.3\textwidth]{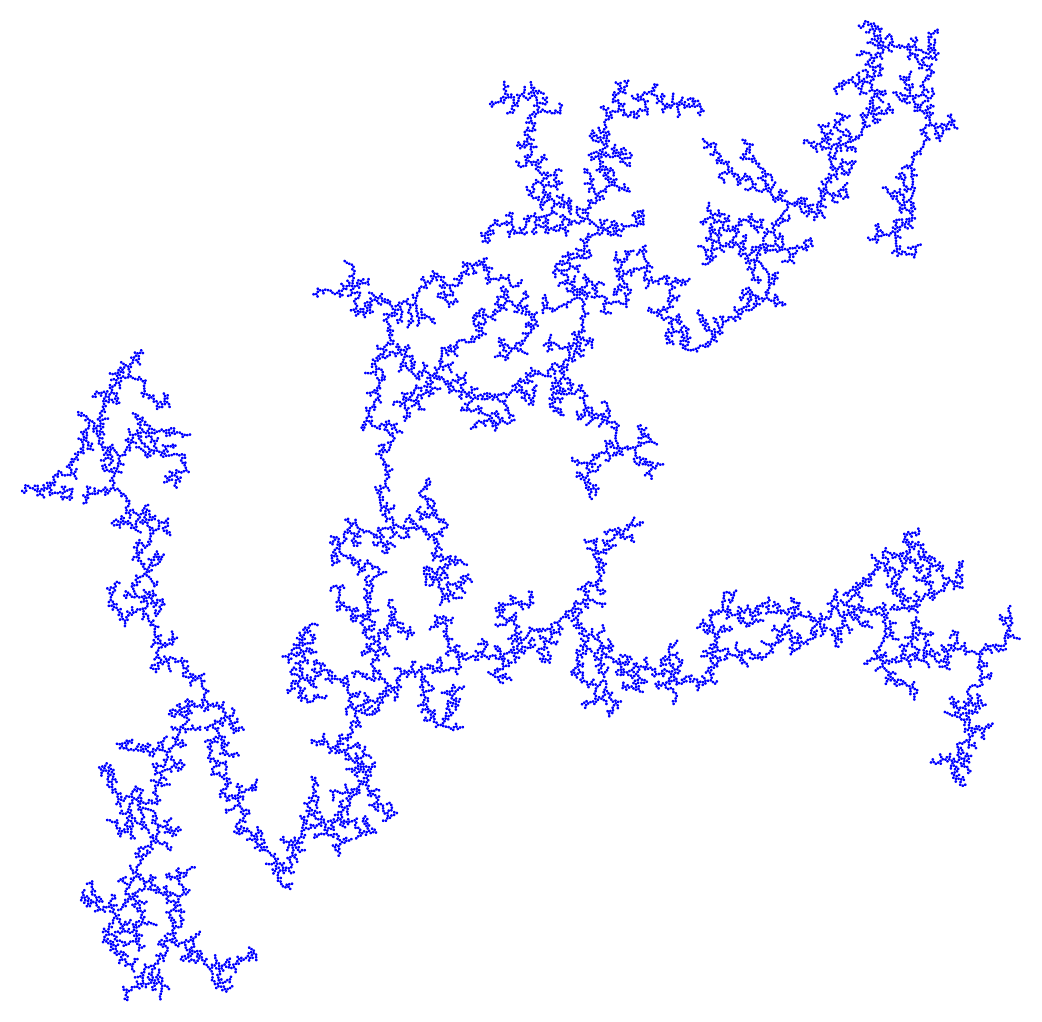}}
\subfigure[Brownian tree]{\label{fig:DLA}\includegraphics[width=.3\textwidth]{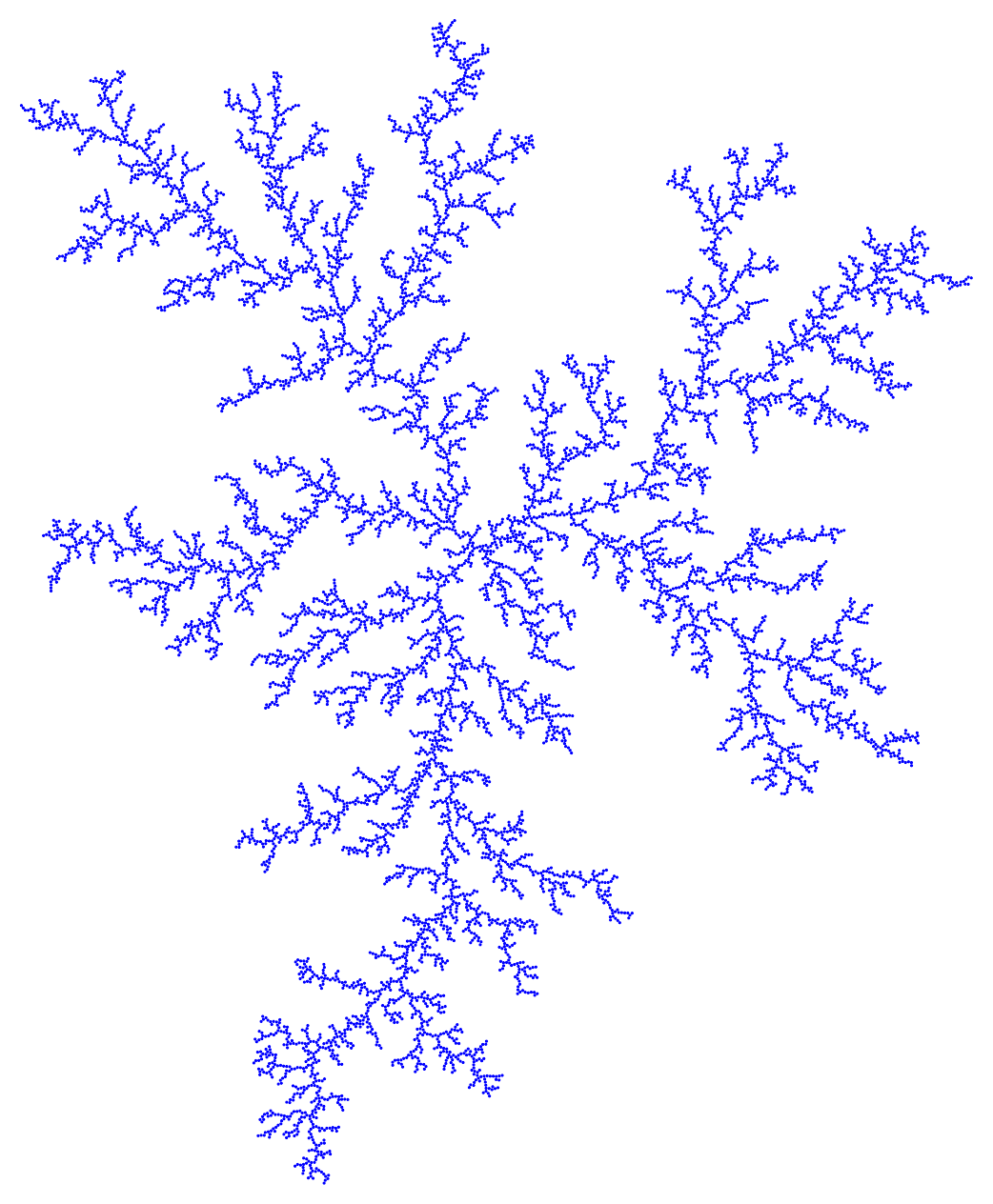}}
\subfigure[Self-avoiding walk]{\label{fig:SAW}\includegraphics[width=.3\textwidth]{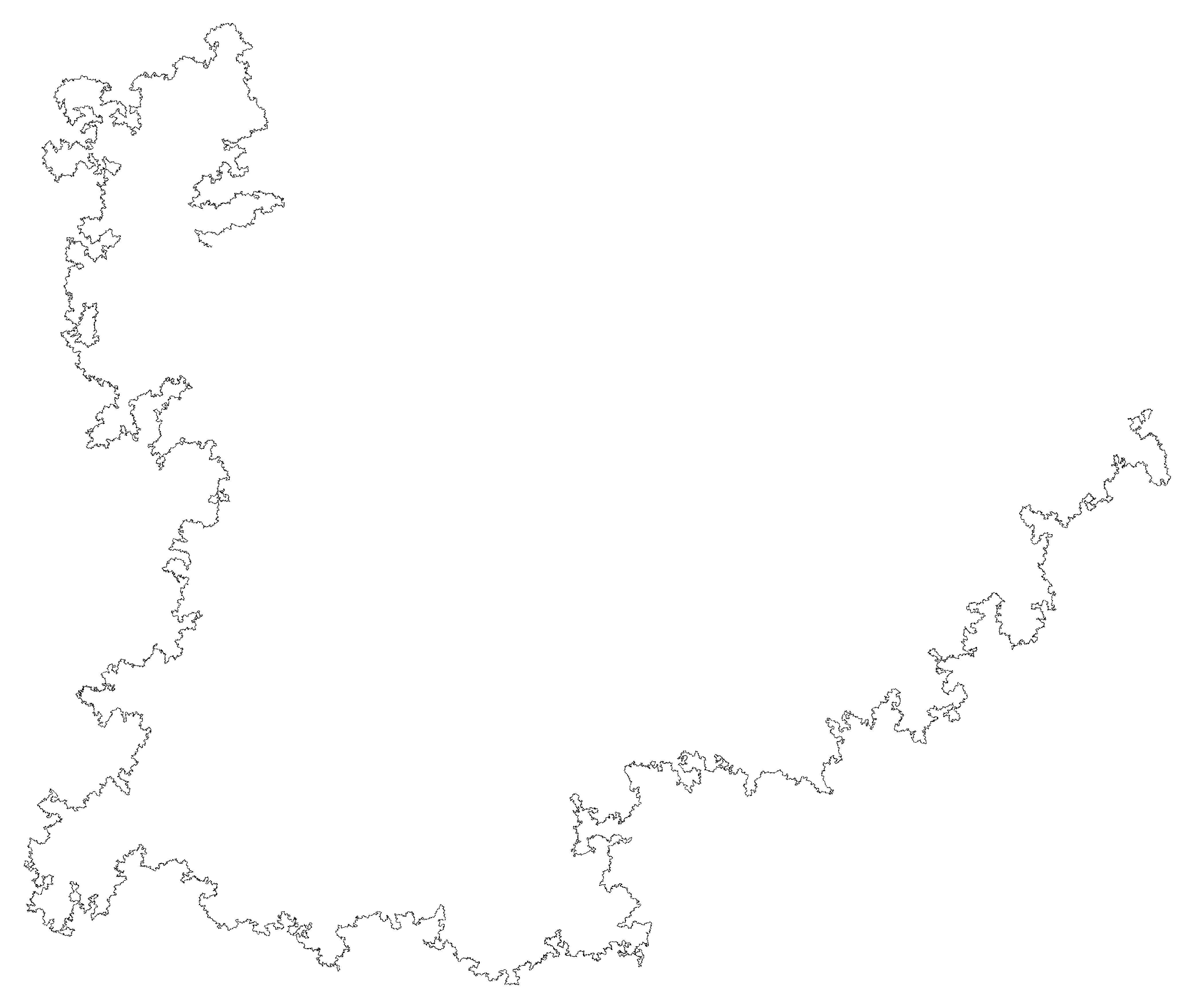}}
\caption{Three examples of polymers constructed by different processes}
\label{fig:polymers}
\end{figure}

We consider geometric structures $S$ in an ambient $m$-dimensional Euclidean.  Many such structures of physical interest have great local complexity, with features on a wide range of length scales. Typically, these structures are generated by a probabilistic process.    Fig.~\ref{fig:polymers} shows the dimension $m=2$ versions of the examples we will consider in detail in this paper.

The dimension $m=3$ versions of these structures are well known models for interesting physical systems.  For example, a branched polymer is a model for a large branched polymer molecule, and a self-avoiding (random) walk is a model for long chain polymer.  A Brownian tree is a  model for soot deposition.  These structures are also represent universality classes of interest to mathematicians.  For example, $m=2$ self-avoiding walk is conjectured to be SLE~{8/3} (where SLE means Schramm Loewner evolution).

In this paper, we propose a mathematical measure of shape suitable for the study of such complicated structures $S$.

\begin{figure}[htp]
\centering
\includegraphics[width=.75\textwidth]{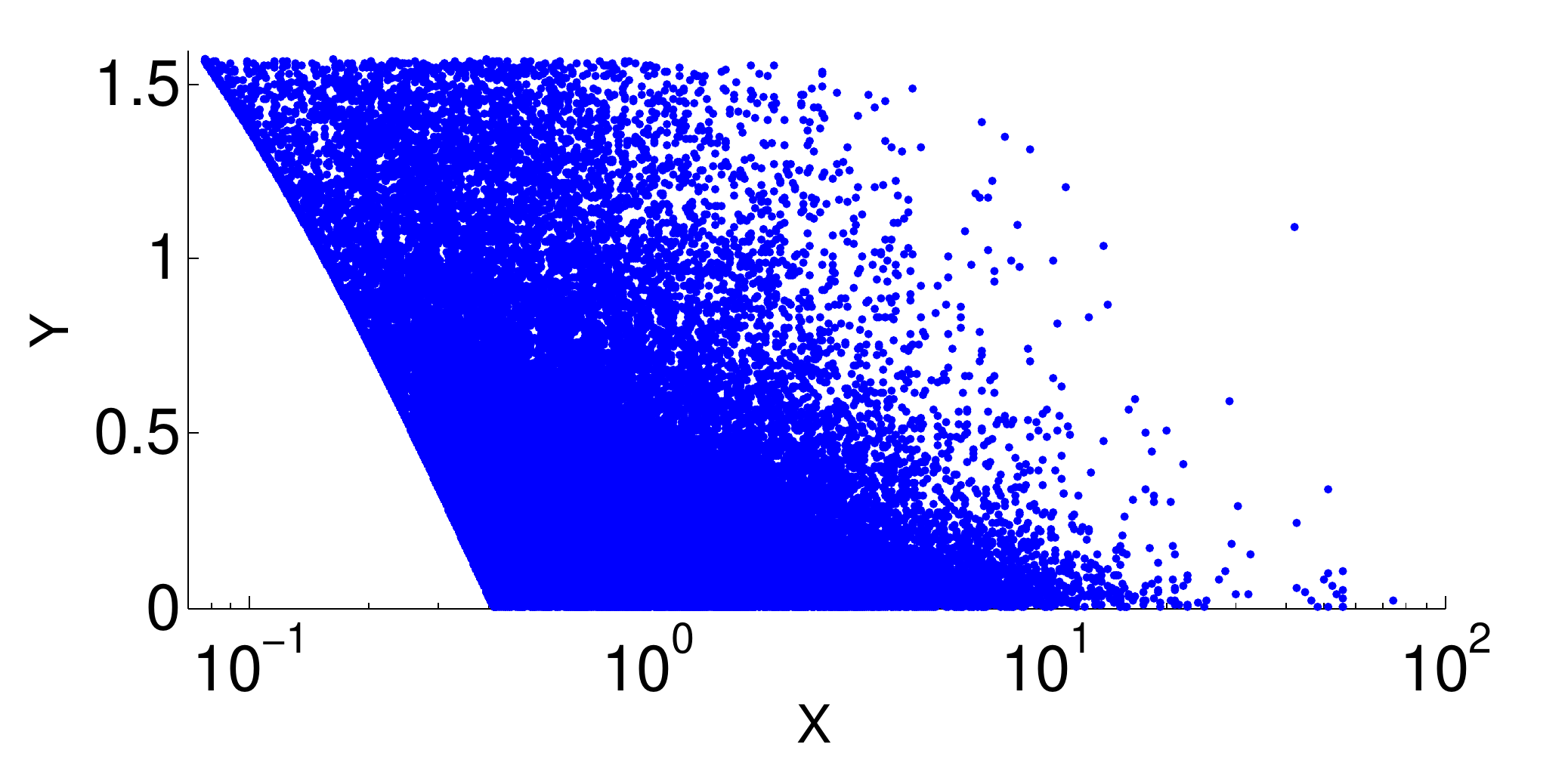}
\caption{Scatter plot of feature points for a branched polymer, $m=2$, $i=1$.}
\label{fig:scatter}
\end{figure}

\newcommand{\PH}{P\!.H. }
\subsection{Distributions of \PH  Points.}  \label{PH Points}We fix an integer $i$, $0 \leq i < m$, which we call the {\em degree}.   Then the geometric structure $S$ gives rise to a scatter plot of points in the plane, as in Fig.~\ref{fig:scatter}.\footnote{The sharp cutoff on the left of the scatter plot is due to the geometry of small configurations of balls, and is not of interest to us.} We  call these points {\em \PH points}, where \PH stands for {\em persistent homology}.  The reason for the terminology will be described later on.  Each \PH point represents a particular kind of $i$-dimensional ``feature'' of $S.$
If the dimension $m$ is $2$ as in Fig.~\ref{fig:polymers} and the degree $i$ is $1$, this feature is intuitively a ``gulf'' \label{gulf} largely surrounded by $S,$ as can be seen in Fig.~\ref{fig:circle}.  The $x$-coordinate of the feature point representing a given gulf
measures its size, and the $y$-coordinate measures what we call its
{\em aspect}: how completely the gulf is surrounded by $S$. Visibly in Fig.~\ref{fig:polymers}, there are many such gulfs of different sizes and aspects, which is reflected in the many \PH points in Fig.~\ref{fig:scatter}.

The distribution of the \PH points in the plane is the measure of shape of $S$ that we study in this paper.  We can define interesting other geometric quantities and properties of $S$ in terms of this  distribution.

\subsection{\PH Dimension and \PH Self-Similarity.} When studying fractal sets $S$ such as the sets mentioned in \S\ref{SETS} above, the notions of Hausdorff dimension and self-similarity are traditionally useful.  These notions have analogues defined in terms of the distribution of \PH points.

We define the {\em \PH dimension} of $S$ by the rate of growth of the density  of \PH points with a given $x$ coordinate as $x$ decreases to $0$ (see \S \ref{Persistent Homology Dimension}).  The \PH dimension is a new kind of fractal dimension of $S$. A  computation shows that in some cases, the \PH dimension agrees with the Hausdorff dimension of $S$.

We say that $S$ is {\em \PH self-similar} if the distribution of points in the plane obtained by rescaling the $x$ coordinate is a scalar multiple of the original distribution of points (for small enough $x$). Ordinary self-similarity is related to Hausdorff dimension. Similarly. \PH self-similarity determines \PH dimension.

\subsection{The Topology of Neighborhoods of $S$.}  The \PH points of $S$ are defined and computed using the topology of neighborhoods of $S$ in the ambient Euclidean space.  For any positive number $\epsilon$, let  $S_\epsilon$ be the  $\epsilon$-neighborhood of $S$, i.e. $S_\epsilon$ is the enlargement of $S$ that includes all points of distance at most $\epsilon$ from $S$.   As  $\epsilon$ grows, the topology of $S_\epsilon$ changes.  We use these changes to define the \PH points.

For example,  if $m=2$ and $i=1$, each \PH point corresponds to a ``gulf'' in $S$.  As $\epsilon$ grows, at some time when $\epsilon=b$, a topological change in $S_\epsilon$ will take place: the mouth of the ``gulf'' fills in so that it becomes a ``lake'' fully enclosed by $S_\epsilon$.  At a later time when $\epsilon=d$, the ``lake'' is filled in.  The number $b$ is the value of $\epsilon$ when the lake is born and $d$ is the value when the lake dies.  The coordinates $x$ and $y$ of the \PH point are defined by the formulas $x=(d+b)/2$ and $y=\arcsec(d/b)$.

The  general definition of the \PH points uses persistent homology.  The ``lake'' in the example above is replaced by a generator of the $i$-dimensional persistent homology of the family $S_\epsilon$ which is born at $\epsilon=b$ and dies at $\epsilon = d$.  The coordinates of the \PH point are determined by the same formulas (see \S \ref{Constructing PH Points with Homology}).  Since many readers will not be at home with homology, we give alternative interpretations of $b$ and $d$ in terms of what kinds of shapes $S$ can ``trap'').

\subsection{Probabilistic Distributions of \PH Points.} When the structure $S$ is generated by a probabilistic process as in the three examples in \S\ref{SETS}, the \PH points are distributed according to a density function $f(x,y)$ determined by that probabilistic process.   Then the \PH dimension of $S$ is given by the growth rate of $f(x,y)$ as $x$ decreases to zero.

We define the process to be {\em \PH statistically self-similar} if the function $f(x,y)$ has the form
\begin{equation*}f(x,y) \approx x^{-(d+1)}\cdot g(y)\end{equation*}
for small enough $x$.  The number $d$ will be the \PH dimension.

\subsection{Computer Calculations.} In most cases, the density function $f$ cannot be calculated in closed form.  However, it can be estimated by creating some sample structures $S$, and then calculating the scatter plot of \PH points for each structure $S$.  This is what we do in this paper.

 We do detailed computer computations for the three geometric structures of Fig.~\ref{fig:polymers}.  In addition to their physical interest, we chose these structures because recently developed  algorithms make computer construction of these structures  tractable.
Also, these structures are topologically trivial, so the fact that we can obtain useful information about them with topological methods is perhaps surprising.  The definition of these structures, together with a description of the algorithms to compute them is in \S\ref{Three Geometric Structures}, and information about the algorithms used to computer their \PH points is in \S\ref{Computational Methods}.

\subsection{Computational Results.} The results of our computations are in  \S\ref{Computational Results}.  They show that
\begin{itemize}\item In cases where the structure $S$ has is known or believed to be statistically self-similar with Hausdorff dimension $d$, we find that it is also \PH statistically self-similar with \PH dimension $d$ (for all degrees $i$).
\item In cases where the structure $S$ fails to be statistically self-similar, it also fails to be \PH statistically self-similar.
\item Interesting differences in the statistical shape of structures $S$ that are unrelated to fractal dimensions can be detected by the $y$ (or aspect) dependency of $f(x,y)$.
\end{itemize}
We conclude that the study of the \PH points of a complex geometric structure provides insight in these cases, and we hope that it will have broad applications in the future.

\newcommand{\reals}{\mathbb{R}}

\section{Constructing the \PH Points.}\label{Constructing PH Points with Homology}  We consider geometric structures $S$ that are subsets of $m$-dimensional Euclidean space $\mathbb{R}^m$. For a number $\epsilon \geq 0$, let $S_\epsilon$ denote the $\epsilon$-neighborhood of $S$, i.e. the set of all points in $\reals^m$ that are of distance at most $\epsilon$ from some point of $S$.  It is useful to think of $\epsilon$ as time, and to think of the system of neighborhoods $S_\epsilon$ as growing over time.

\subsection{The Definition Using Persistent Homology.} Readers unfamiliar with  homology should skip this formal definition, and go directly to \S\ref{An Example} and \S\ref{Alexander Duality} where the idea is described more intuitively.  For an integer $i$ with $0\leq i <n$, $\tilde H_i(S_\epsilon,\reals)$ denotes the degree $i$ reduced homology group with real coefficients, and is a real vector space measuring the topology of the $\epsilon$-neighborhood $S_\epsilon$.  As $\epsilon$ varies, we get a family $\mathcal{F}$ of vector spaces, one for each value of $\epsilon$.  The fundamental result of persistent homology~\cite{zcarlsson} is that the family $\mathcal{F}$ is a sum of extremely simple families $\mathcal{F}_j$ each of which represents a one dimensional vector space that is born at time $\epsilon=b_j$ and dies at time $\epsilon=d_j$.\footnote{This classification holds for any nice subset of Euclidean space and, in particular, for the geometric structures we consider here.}  Our scatter plot of \PH points has one \PH point for each simple family  $\mathcal{F}_j$ are our {\em features} with coordinates  $(x,y)$ where
\begin{equation}\label{coordinates} \text{the {\em size}}\;\;x = \left(\frac{b_j+d_j}{2}\right) \;\;\;\;\;\;\;\;\text{and the {\em aspect}}\;\;y=\arcsec\left(\frac{d_j}{b_j}\right)\end{equation}

If $S$ is connected (as in the three examples in Fig.~\ref{fig:polymers}), then $\tilde H_0(S_\epsilon)$ vanishes for all $\epsilon$, so there are no  \PH points for $i=0$.  If $S$ has diameter $\delta$, i.e. the distance between any two points of $S$ is $\leq \delta$, then the $x$ coordinate of any \PH point is $\leq \delta$ because $\tilde H_i(S_\epsilon)$ vanishes for $\epsilon >\delta$. \label{delta}

\subsection{An Example}\label{An Example}
\begin{figure}[htb]
\centering
\subfigure[$\epsilon=0$]{\label{fig:here1}\includegraphics[width=.195\textwidth]{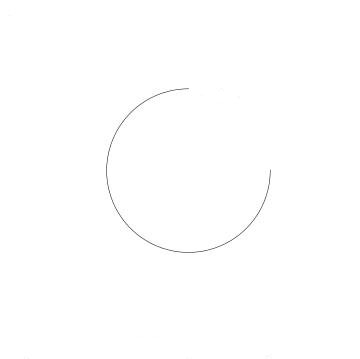}}
\subfigure[$0<\epsilon<\cos(\theta)$]{\label{fig:here2}\includegraphics[width=.185\textwidth]{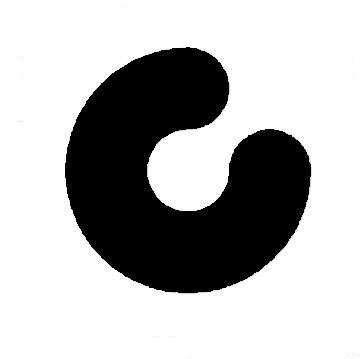}}
\subfigure[$\epsilon=\cos(\theta)$]{\label{fig:here3}\includegraphics[width=.185\textwidth]{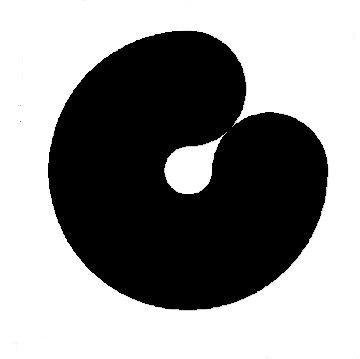}}
\subfigure[$\cos(\theta)<\epsilon<1$]{\label{fig:here4}\includegraphics[width=.185\textwidth]{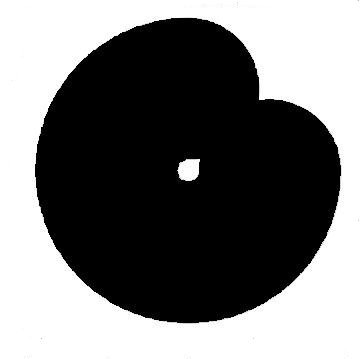}}
\subfigure[$\epsilon=1$]{\label{fig:here5}\includegraphics[width=.185\textwidth]{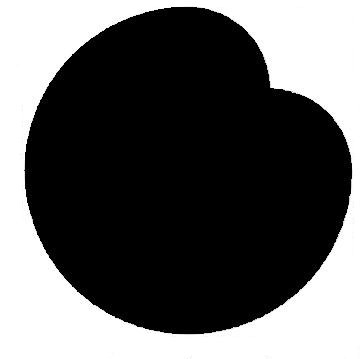}}
\caption{Neighborhoods of a circle with an arc removed (r=1)}
\label{fig:circle}
\end{figure}

Suppose that $S$ is an arc of a circle of radius $r$ in the plane.  Suppose that the angle of the arc is $\pi + 2\theta$ (so its length is $r(\pi + 2\theta)$).  Now we can draw the $\epsilon$-neighborhoods $S_\epsilon$ for various $\epsilon,$ as in Figure~\ref{fig:circle}.

We see that the topology of $S_\epsilon$ doesn't change until $\epsilon=r\cos(\theta)$, at which time it develops a hole.  This hole persists until $\epsilon=r$, when the hole if filled in.  The hole gives rise to an element of $H_1(S_\epsilon)$, which is one of our ``features'' $\mathcal{F}_j$ for $i=1$.  It is born at time $\epsilon=r\cos(\theta) = b_j$ and it dies at time $\epsilon=r=d_j$.  It gives rise to a \PH point $(x,y)$, where
$$ x = \left(\frac{b_j+d_j}{2}\right) =  \left(\frac{r+r\cos\theta}{2}\right) \;\;\;\;\;\;\;\;\text{and}\;\;\;\;\;\;\;\;\; y=\arcsec\left(\frac{d_j}{b_j}\right)=\theta  $$

As in \S\ref{gulf}, the circular arc is like a gulf.  The \PH point's $x$-coordinate is proportional to $r$, so it is a measure of the size of the gulf.  Its $y$-coordinate or aspect is $\theta$ which measures the or how completely the gulf is surrounded by $S$.  For example, when $\theta$ takes it maximum value of $\pi/2$, the gulf is completely surrounded. This example easily generalizes to the two-sphere in $\R^3$ with a disk of angular size $\alpha$ removed.

\subsection{Interpretation by Objects in Cages.}\label{Alexander Duality}  There is an alternative interpretation of \PH points in terms of bodies of various types trapped or caged by $S$.

Consider a solid structure $S$ fixed in  Euclidean $3$-space and a mobile spherical ball $B$ of variable radius $\epsilon$ that cannot pass through $S$.  The ball could be trapped or \em{caged} \rm by a gulf or hole of $S.$  If the ball expands, at some point $\epsilon = d_j$ the ball will no longer fit into its cage, so the cage dies for $\epsilon$-balls as $\epsilon$ grows beyond $d_j$.  If  $\epsilon$ decreases so the ball shrinks, at some point $\epsilon=b_j$ it can escape it cage, either into the exterior or into another cage, so the cage is born only when $\epsilon$ grows beyond $b_j$.  The numbers $b_j$ and $d_j$  coincide with the birth and death values of a persistent homology class for degree $i=2$ (and $m=3$).  The reason is Alexander duality.  As always, the \PH point $(x,y)$ corresponding to this cage is computed by equation (\ref{coordinates}).

In dimension $m=2$, for degree $i=1$ the features for $S$ can be similarly interpreted as cages trapping a disk of radius $\epsilon$.  When $\epsilon<b_j$ the disk can escape the cage, and when $\epsilon>d_j$ the disk won't fit in the cage.

For any $i$ and $m$, Alexander Duality yields an interpretation of $H_i(S_\epsilon)$ in terms of $\epsilon$-neighborhoods of mobile surfaces of dimension $m-i-1$ trapped by $S$.  For example, if $i=1$ and $m=3$, instead of trapping balls or radius $\epsilon$, we are trapping snakes of radius $\epsilon$ (i.e. the $\epsilon$-neighborhood of a curve).  We note that if $S$ is a physical polymer, its ability to trap particles of size $\epsilon$ or tubes of radius $\epsilon$ is of physical interest.

\section{Dimension and Self-Similarity.}\label{Persistent Homology Dimension}

In this section, we introduce three notions that are derived from the scatter plot of \PH points: \PH dimension and \PH self-similarity, and \PH statistical self-similarity.  These ideas are particularly useful when  $S$ is fractal.  They are \PH analogues of  Hausdorff dimension and ordinary self-similarity.

\subsection{\PH Dimension and \PH Self-Similarity.} \label{F} Choose a degree $i$  and a subinterval $I$ of the interval from $0$ to $\pi/2$. (If no interval is specified, we take $I$ to be the whole interval from $0$ to $\pi/2$.) Let $F(x)$ be the number of \PH points $(x_j,y_j)$ of $S$ such that $x_j \geq x$ and $y_j$ is in the interval $I$.\footnote{For the applications discussed in this paper, we always take $I$ to be the whole interval from $0$ to $\pi/2.$ However, we believe that the choice of an interval $I$ will be important in other applications.} Since $F$ is essentially an integral from the right with respect to $x$,  $F$ is a nowhere increasing function.

\subsubsection{Definition of \PH Dimension.}  \label{persistent homology dimension} We say that $S$ has { \PH dimension} $d$ (for the data $i$ and $I$) if
\begin{align} \text{for all} \;\;\;c<d,\;\; \lim_{x\rightarrow0} x^cF(x) &=\infty \;\;\;\;\;\text{and}\\\text{for all} \;\;\;c>d,\;\; \lim_{x\rightarrow0} x^cF(x) &= 0\end{align}
i.e. if $F(x)$ grows roughly like $x^{-d}$ as $x$ decreases to $0$.
\newcommand{\rn}{\ell}
\subsubsection{Definition of \PH Self-Similarity.} Let $\rho$ be a number, $0<\rho<1$, and let $\rn$ be an integer, $\rn>0$.  We say that $S$ is {\em \PH self-similar} with similarity ratio $\rho$ and replication number $\rn$ if for small enough $x$, the number of \PH points with coordinates $(\rho x,y)$ is $\rn$ times the number of \PH points with coordinates $(x,y)$.  (To rule out trivial cases, we assume $S$ has \PH points for arbitrarily small $x$ coordinate.)

\subsection{The Relation between \PH Dimension and \PH Self-Similarity.} The following proposition is an analogue of the relationship between self-similarity and Hausdorff dimension.

\subsubsection{Proposition.}  If $S$ is \PH self-similar  with similarity ratio $\rho$ and replication number $\rn$, then its \PH dimension is
$$d=\frac{\log( \rn)}{\log( \rho^{-1})}=-\frac{\log \rn}{\log \rho}$$
(To rule out trivial cases, we assume $S$ has \PH points for arbitrarily small $x$ coordinate, whose $y$ coordinate is in $I$.)

\subsubsection{Proof.} Suppose that $x\leq c$ is ``small enough x'' in the definition of \PH self-similar. It is enough to show that for $x\leq c$, $A_1+B_1x^{-d} \leq F(x)\leq A_2+B_2x^{-d}$ for some $A_i,B_i$ and for $d=-\log \rn/\log\rho$.

For positive integers $k$,
$$F(\rho^k c) = F(c) + \rn + \rn^2+\cdots+\rn^k = A+B(\rho^k)^{-d}$$ where $A=F(c)+\rn/(1-\rn)$, and $B=\rn/(1-\rn)$.
Now suppose $\rho^{k+1}c\leq x\leq \rho^k c$.  Then
$$F(x)\geq F(\rho^k c)=A+B(\rho^k c)^{-d}\geq A+B(x/\rho)^{-d} = A+B\rho^{d}x^{-d}\equiv A_1+B_1x^{-d}$$
$$F(x)\leq F(\rho^{k+1} c)=A+B(\rho^{k+1} c)^{-d}\leq A+B(x \rho)^{-d} = A+B\rho^{-d}x^{-d}\equiv A_2+B_2x^{-d}$$

\subsection{Relation between \PH dimension and Hausdorff dimension.} The formula $d=-\log \rn / \log \rho$ shows that \PH dimension bears the same relation to \PH self-similarity that Hausdorff dimension bears to ordinary self-similarity.  Nevertheless, \PH dimension and Hausdorff dimension in general measure different things.  Here are some examples:
\subsubsection{Hausdorff dimension can exceed \PH dimension.}  $d$-dimensional Euclidean space has Hausdorff dimension $d$ but it has \PH dimension $0$ (since the reduced homology of Euclidean space vanishes, so $F(x)=0$ for any $i$ and $I$).

This example highlights a difference: \PH dimension always measures complexity, whereas Hausdorff dimension measures some combination of complexity and topological dimension.

\subsubsection{\PH dimension can exceed Hausdorff dimension .} For any $\rho$ and $\rn$, let the set  $S_{\rho,\rn}$ be  the set of real numbers
$$0=x_{1,0},x_{1,1},\ldots,x_{1,\rn}=x_{2,0}, x_{2,1},\ldots,x_{2,\rn^2}=x_{3,0},x_{3,1}\ldots, x_{3,\rn^3}=x_{4,0}, \ldots$$ where $x_{k,j+1}-x_{k,j} = 2\rho^k$.  The set $S_{\rho,\rn}$ has \PH points only for $i=0$, where the points are all of the form $(\rho^k,\pi/2)$ which occurs $\rn^k$ times for each  $k=1,2,\ldots$ (representing a degree 0 persistent homology class that is born at $\epsilon=0$ and dies at $\epsilon=2\rho^k$).  So it is \PH self-similar  with similarity ratio $\rho$ and replication number $\rn$.  Therefore, it has \PH dimension $-\log \rn / \log \rho$ if $i=0$ and $I$ contains $\pi/2$.  However,  it has Hausdorff dimension $0$ and it is not self-similar in the usual sense.

\subsection{An Example that is Both Self-Similar and \PH Self-Similar.}

We present a fractal $S$ in the plane that is interesting because  its \PH points can be calculated exactly, and its Hausdorff dimension agrees with its \PH dimension.

\subsubsection{Construction of $S$.} 
\begin{figure}[htp]
\centering
\includegraphics[width=\textwidth]{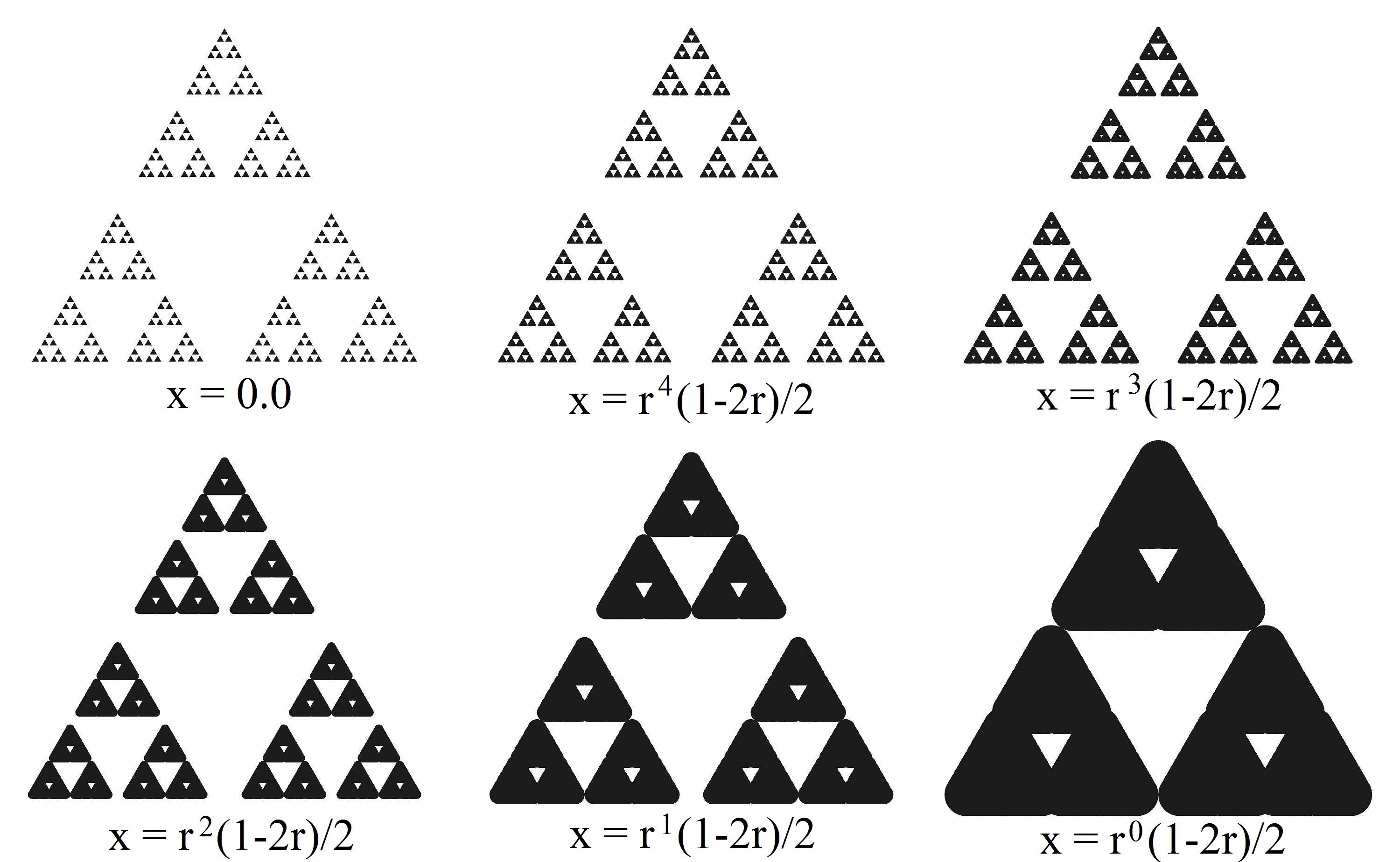}
\caption{A self-similar fractal and its $\epsilon$-neighborhoods}
\label{fig:triangles}
\end{figure}

The fractal $S$ is depicted with several its $\epsilon$-neighborhoods in Figure ~\ref{fig:triangles}. It is a generalized version of the Sierpinski triangle. To construct it, let $S^1$ be an equilateral triangle with sides of length one and vertices $\paren{0,0},\paren{1,0},$ and $\paren{1/2,\sqrt{3}/2}.$ Let $f_1,f_2,$ and $f_3$ be contractions of the plane with shared factor $\rho\leq1/2$ that fix $\paren{0,0},$ $\paren{1,0}$, and $\paren{1/2,\sqrt{3}/2},$ respectively. Let $S^k=f_1\paren{S^{k-1}}\cup f_2\paren{S^{k-1}}\cup f_3\paren{S^{k-1}}$ for $k>1.$ Then  $S$ is defined to be $\cap_{k}S^k.$ If $\rho=1/2,$ $S$ is the Sierpinski triangle.

\subsubsection{Proposition.} $S$ is \PH self-similar with similarity ratio $\rho$ and replication number $3$.  Therefore its \PH dimension is $-\log 3/\log\rho$ (provided that either  $i=1$ and the interval $I$ contains the number $\tilde
y$ defined below, or else  $i=0$ and the interval $I$ contains the number $\pi/2$). Furthermore, it is self-similar in the usual sense with similarity ratio $\rho$ and replication number $3$, so its Hausdorff dimension is also $-\log 3/\log\rho$.

\subsubsection{ Proof.} We calculate all of the \PH points of $S$ exactly.

First, consider $i=1$.  For every integer $k\geq 0$, there are $3^k$ independent cycles in $H_1(S_\epsilon)$ which are born when  $\epsilon={\rho^{k}\paren{1/2-\rho}}\equiv b_k$.  These are the values of $\epsilon$ illustrated in Figure~\ref {fig:triangles}, where the cycles are visible as newly formed holes in $S_\epsilon$.  These cycles die when $\epsilon=\rho^k\sqrt{1/3-\rho+\rho^2-\rho^3+\rho^4}\equiv d_k.$ Therefore $S$ has $3^k$ \PH points $(x_k,y_k)$ where  $x_k=(b_k+d_k)/2=\rho^k(b_0+d_0)/2=\rho^kx_0$ and $y_k=\text{arcsec}\paren{d_k/b_k}=\text{arcsec}(d_0/b_0)\equiv \tilde y$ is independent of $k$.

Next consider $i=0$.   For every integer $k\geq 0$, there are $2\cdot3^{k}$ cycles that die when  $\epsilon={\rho^{k}\paren{1/2-\rho}}\equiv d_k$. A canonical basis of $H_0(S_\epsilon)$ is indexed by the connected components of $S_\epsilon$.  So we can see the homology classes dying since  connected components are fused together at these values of $\epsilon$, Fig. \ref {fig:triangles}.   Every homology class is born at $\epsilon=0$, so we take $b_k=0$ for all $k$.    Therefore $S$ has $2\cdot3^{k}$ \PH points $(x_k,y_k)$ where where $x_k=(b_k+d_k)/2=d_k/2=\rho^k(d_0)/2=\rho^kx_0$ and $y_k=\text{arcsec}(d_k/b_k)=\text{arcsec}(\infty)=\pi/2 $, again independent of $k$.

\subsection{Probabilistically Determined Structures $S$.}

Unlike the example above, the physically interesting structures from \S\ref{SETS} we consider are determined by a probabilistic process.  We state versions of self-similarity and dimension adapted to this situation.

\subsubsection{The \PH Density Function.}  \label{distribution} The \PH points are drawn from a distribution with a {\em \PH density function} $f(x,y)$ which measures the expected number of \PH points per unit area at $(x,y)$.  The density function $f$ depends only on the probabilistic process giving rise to $S$, unlike the individual \PH points which depend on a particular realization $S$ of that process.
In this case, the expected value $\mathbb{E}\left[F(z)\right]$ of the function $F(z)$ will be
$$ \mathbb{E}\left[F(z)\right]= \int_{y\; \text{in}\; I} \int_{x=z}^\infty f(x,y)\,\text{d}x\,\text{d}y
$$

\subsubsection{\PH Statistically Self-Similar Processes.} \label{hypothesis} We call the process giving rise to $S$ \PH statistically self-similar with self-similarity dimension $d$ if the \PH points in the plane are sampled from a spatial Poisson process with  \PH density function of the form
\begin{equation*}f(x,y) \approx x^{-(d+1)}\cdot g(y)\end{equation*}
for some function $g(y)$, at least for $x$ small compared to the diameter $\delta$ of $S$ (at which point $f(x,y)$ will vanish).  Here $g(y)$, which is a function of the aspect $y$ alone, encodes additional information about the shape of $S$ unrelated to the dimension.
In this case,
\begin{equation*}\label{hypothesis F} \mathbb{E}\left[F(z)\right] \approx z^{-d}\end{equation*} up to a multiplicative constant, independent of the interval $I$, and  $d$ will be the \PH dimension almost surely.

\subsubsection{Hypothesis.} We hypothesize that many physically interesting statistical structures $S$ that are statistically self-similar in the usual sense with similarity dimension $d$  are also \PH statistically self-similar with self-similarity dimension $d$.  We give numerical evidence for this hypothesis in \S\ref{Computational Results} in the case of 2D self-avoiding walks and 3D branched polymers by  calculating the \PH points of individual realizations $S$ of that process (averaging to improve accuracy) and estimating  $F$ from this data.

\section{Three Geometric Structures}\label{Three Geometric Structures}

An $m$-dimensional polymer of order $n$ is a connected subset of $\mathbb{R}^m$ that is the union of $n$ distinct $m$-balls of the same radius $r$.  We assume that the balls are non-overlapping, i.e. two $m$-balls are either disjoint, or else tangent.  In this paper, we consider the objects drawn from three different probability distributions on the space of polymers. These objects are called branched polymers, Brownian trees, and self-avoiding walks.

We have chosen these three structures to study for several reasons.  First, they are all models for physical processes of interest, and they have been studied by mathematical physicists as representatives of interesting universality classes.  Second, in each case the union of the balls is itself contractible with probability $1$, so the fact that topology can yield interesting information is perhaps surprising.  Third, there has been a theoretical breakthrough leading to an algorithm to efficiently compute random representatives drawn from the probability distribution in each case.

\subsection{Branched Polymers.}  Branched polymers are drawn from a distribution that can be viewed as the uniform distribution on the space of polymers. Given a  polymer composed of $n$ unit balls $\set{B_1,\ldots,B_n}$ each of radius $r$, let $x_1,\ldots,x_n$ be their centers.  We define the tree-type of the polymer to be the tree with vertices $\set{1,\ldots,n}$ and edges $\set{i,j}$ such that $B_i$ and $B_j$ are tangent.  There is a measure $\Omega_T$ on the space of branched polymers whose tree-type is $T:$
$${\Omega_T= \prod_{{ij\in T}\atop{i<j}} \Omega\paren{x_i-x_j}}$$
where $\Omega\paren{x_i-x_j}$ is the area form of the unit sphere of vectors $x_i-x_j$.  The probability measure on the space of all branched polymers is defined as the sum of this measure over each tree type, normalized to have total volume one.

We used an algorithm created by Kenyon and Winkler to inductively construct branched polymers in two and three dimensions~\cite{kenyon}. A two-dimensional branched polymer of order 10,000 is depicted in Fig.~\ref{fig:BP}.

\subsection{Brownian Trees.} Another distribution on the space of polymers samples objects called Brownian trees. They are built inductively -- given a Brownian tree with $n-1$ balls, another ball is placed on a circle at some distance away from it. It then moves via Brownian motion until it reaches the existing structure and sticks. This process is called diffusion limited aggregation.  A Brownian tree composed of 10000 balls is depicted in Figure~\ref{fig:DLA}.

To create Brownian trees, we used a program written by Mark J. Stock~\cite{stock}.

\subsection{Self-Avoiding Walks.} The third distribution we consider is that of self-avoiding (random) walks on a square lattice, which are a model for linear polymers. A self-avoiding walk is a subset of the edges of a lattice that forms a non-intersecting path. They are given the uniform distribution: if there are $\ell$ self-avoiding walks of length $n,$ each is given weight $1/\ell.$ We realize self-avoiding walks as polymers by viewing each edge as three touching balls. An example of a self-avoiding walk composed of 500,000 lattice edges is shown in figure~\ref{fig:SAW}.

We used the fast-pivot algorithm created and implemented by Tom Kennedy to generate self-avoiding walks~\cite{pivot}~\cite{kennedy}.

\subsection{The Limit.}\label{The Limit} The polymers we consider are composed of finite unions of balls of a fixed radius and therefore are of finite complexity.  If we take the limit as the size of the balls goes to $0$ and the number of balls increases so the expected diameter of the polymer remains fixed, we obtain objects of mathematical interest. One wants dimension estimates of these limiting objects and, in particular, we would like to compute their persistent homology dimensions (if they exist).

One expects that the large scale properties of the polymers  remain unchanged once the ball radius $r$ is small enough, and that these reflect the properties of the limiting object.  In particular, we expect the distribution of \PH points $(x,y)$ to be independent of $r$ once the size coordinate $x$ is large with respect to $r$.

\section{Computing Persistent Homology of Polymers.}\label{Computational Methods}

We calculate the persistent homology of a polymer by first computing filtrations of alpha complexes and then loading this information into the JPlex computational homology software package ~\cite{jplex}.

The alpha complex $X_{\alpha}$ of a set of points $X$ is a simplicial complex homotopy equivalent to the union of the set of balls of radius $\alpha$ centered at those points. It was developed by Herbert Edelsbrunner and others~\cite{edelsbrunner1}~\cite{edelsbrunner2}~\cite{edelsbrunner3}.  The alpha complex is a subcomplex of the Delaunay triangulation. To be specific, it consists of all simplices $\sigma$ in the Delaunay triangulation of $X$ such that the vertices of $\sigma$ lie on a sphere of radius $\leq\alpha$ that is empty in $X.$  Refer to~\cite{Harer} for more information concerning the alpha complex and details about the algorithm used to compute it.

Suppose $X$ is a polymer composed of balls centered at points $\set{x_i}$ of radius $r.$  In order to compute the persistent homology of this polymer, we compute a sequence of alpha complexes of $\set{x_i},$  $A_r \hookrightarrow A_{\epsilon_1} \hookrightarrow \ldots \hookrightarrow A_{\epsilon_n},$ where each $\epsilon_i$ is a value at which a new simplex $\paren{\text{or simplices}}$ is added to the alpha complex and $A_{\epsilon_n}$ contains every simplex in the Delaunay triangulation. This computation allows us to assign a filtration value to each simplex in that triangulation -- the lowest value $\epsilon$ such that the simplex appears in the alpha complex $A_{\epsilon}$.

The set of simplices together with their filtration values is then loaded into the JPlex software developed by mathematicians and computer scientists at Stanford University. This program computes the persistent homology of the polymer from this data~\cite{jplex}. JPlex outputs a set of intervals $\paren{b_j,d_j}$ corresponding to the lifetime of a persistent homology class of the polymer during the $\epsilon$-expansion. This data is then processed to create a scatter plot of \PH points.

\section{Computational Results.}\label{Computational Results}

\subsection{Numerical Determination of \PH Statistical Self-Similarity and Estimation of \PH Dimension.} We consider  2D self-avoiding walks, 2D and 3D branched polymers, and 2D and 3D Brownian trees. Using the methods described in the previous section, we compute scatter plots of \PH points for each polymer. Then, we calculate the function $F\paren{x}$ defined in Section~\ref{F} and observe whether or not it appears to follow a power law for some interval on the $x$ axis. If so, we use linear regression (taking the log of both axes so that a power law $F(x)=cx^{-d}$ becomes a linear relation) to estimate the exponent, which will be an estimate of the \PH dimension.

\begin{figure}[htp]
\centering
\subfigure[$F(x)$ for degree $i=1$]{\label{fig:BP3D1}\includegraphics[width=.45\textwidth]{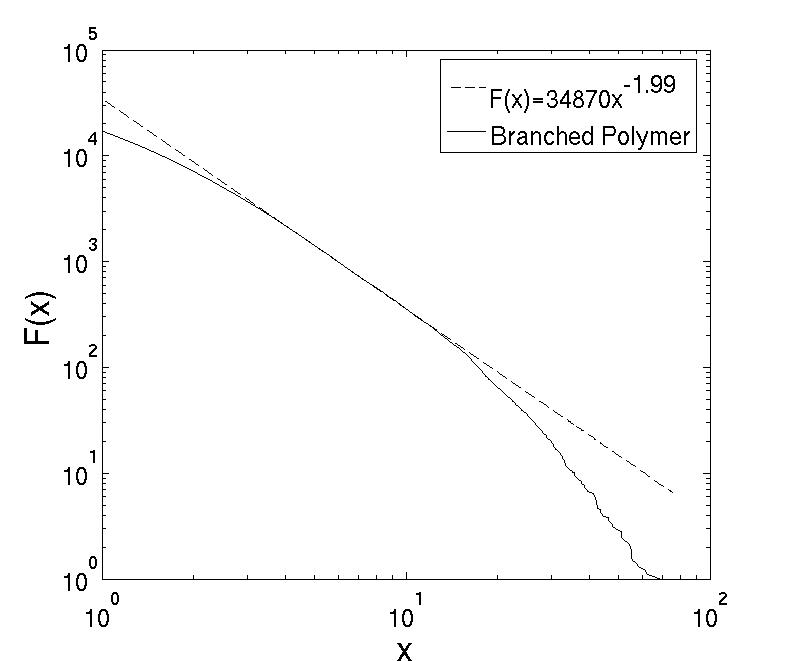}}
\subfigure[$F(x)$ for degree  $i=2$]{\label{fig:BP3D2}\includegraphics[width=.45\textwidth]{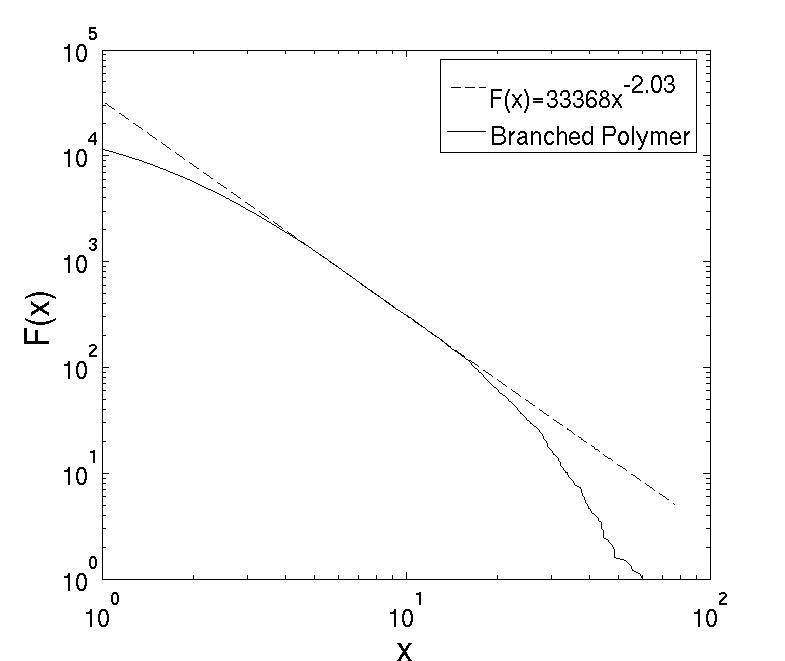}}
\caption{Log-log plots of $F(x)$ for 3D branched polymers}
\label{fig:BPRectanglePlots}
\end{figure}

\begin{figure}[htp]
\centering
\subfigure[2D S.A.W. using $H_1$]{\label{fig:SAWPlot}\includegraphics[width=.45\textwidth]{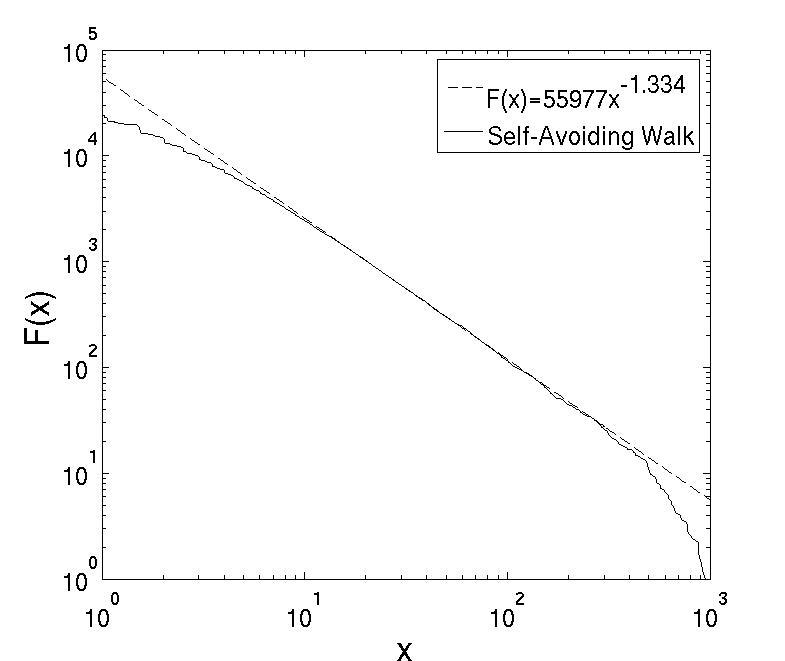}}
\subfigure[2D B.P. using $H_1$]{\label{fig:BP2D}\includegraphics[width=.45\textwidth]{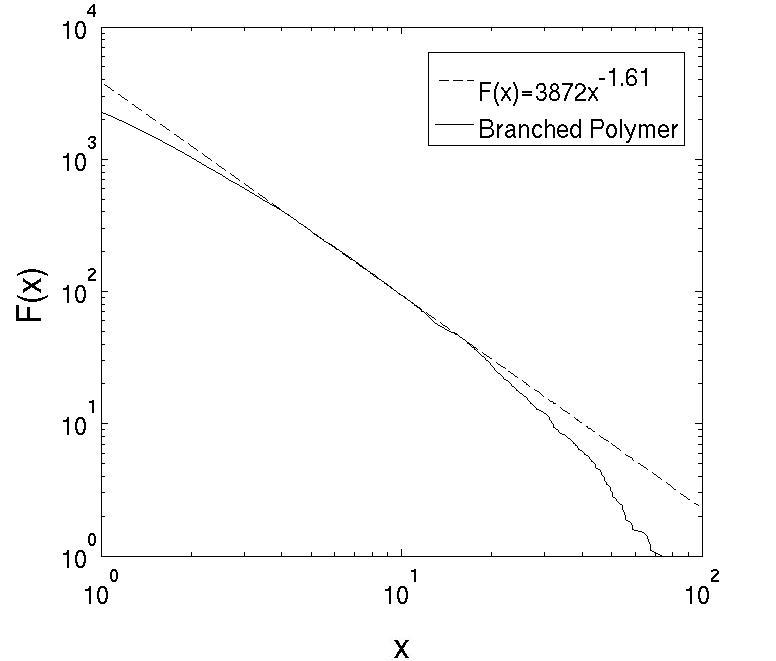}}
\caption{Size/Features plots for 2D self-avoiding walks and 2D Branched Polymers}
\label{fig:BPSAW2DRectanglePlots}
\end{figure}

Even if the probabilistic process we are modeling is \PH statistically self-similar,  we do not expect $F\paren{x}$ to follow an power law for either very high values or very low values of x in our individual simulations. Each simulation has two natural length scales: the diameter $\delta$ of the polymer $S$ itself, and the diameter $r$ of an individual ball in the polymer.  We can only expect that our estimate of $F(x)$ is of physical interest if $x$ is well above $r$ and well below $\delta.$ The function $F(x)$ is identically zero when $x>\delta$ (\S\ref{delta}) so it cannot follow a power law.   Increasing the number $n$ while fixing $r$ would increase the upper end of the range of accuracy.   When $x$ is close to $r$, the features reflect the local geometry of small collections of balls, rather than the limiting object. Note that power law behavior begins at a small multiple of the component ball size in Figures~\ref{fig:BPRectanglePlots} and~\ref{fig:BPSAW2DRectanglePlots} ($r=1$ for branched polymers and Brownian trees, while it is 1/2 for self-avoiding walks which have additional local structure as they are constructed on a lattice).

 Our finite polymers do not have a non-zero \PH  dimension as defined in Section~\ref{persistent homology dimension}, because that takes the limit as $x\rightarrow 0$.

\subsection{Three Dimensional Branched Polymers.}  Plots of $F\paren{x}$ for degree $i=1$ and degree $i=2$ of three-dimensional branched polymers is shown in Figure~\ref{fig:BPRectanglePlots} on log-log paper. The data were gathered from ten branched polymers of order 30,000. For both $i=1$ and $i=2$, the plot appears to be linear from about size $\sim 4$ to size $\sim 15$, so the process appears to be \PH statistically self-similar.    Linear regression yields exponents $-1.99$ and $-2.03,$ respectively. Thus the \PH  dimension of 3D branched polymers is about $2$, which is consistent with the scaling exponent of $1/2$ that Brydges and Imbrie proved for three-dimensional branched polymers in~\cite{brydges}.

\subsection{Two Dimensional Branched Polymers.} Here, our results are inconclusive.  The function $F\paren{x}$ for 2D branched polymers appears to follow a power law with exponent $-1.61$ for a certain range. It was computed from a data set of 10 branched polymers of order $n=10,000$ and it is plotted in figure~\ref{fig:BP2D}. However, the fit is not as convincing as it is for 3D branched polymers and 2D self-avoiding walks.  Furthermore, the rigorous methods of Brydges and Imbrie~\cite{brydges} do not work in two dimensions and do not provide a value for the scaling exponent of 2D branched polymer.

\subsection{Two Dimensional Self-avoiding Walks.}  A similar plot for the  2-dimensional self-avoiding walks and degree $i=1$ is shown in figure~\ref{fig:SAWPlot}. The data were gathered from five self-avoiding walks composed of 500,000 lattice edges each. The plot appears to follow a power law for a very long range: from size $\sim10$ to size $\sim150$, which is consistent with \PH statistical self-similarity.  Least squares regression yields an estimate of $-1.334$ for the exponent and therefore a value of 1.334 for the \PH dimension. This is very close to the conjectured value of $4/3$ for the Hausdorff dimension, which was first hypothesized by Chemist P. Flory and has been proven rigorously for an object believed to be in the same universality class as two-dimensional self-avoiding walks by B. Nienhaus~\cite{nienhuis}.

\begin{figure}[htp]
\centering
\subfigure[$F(x)$ for degree $i=1$]{\label{fig:BT3D1Plot}\includegraphics[width=.45\textwidth]{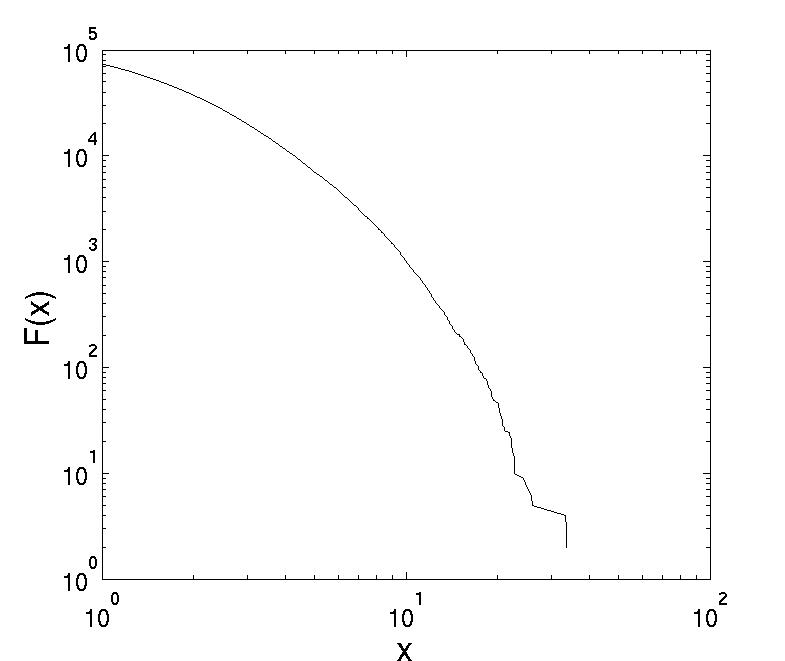}}
\subfigure[$F(x)$ for degree $i=2$]{\label{fig:BT3D2Plot}\includegraphics[width=.45\textwidth]{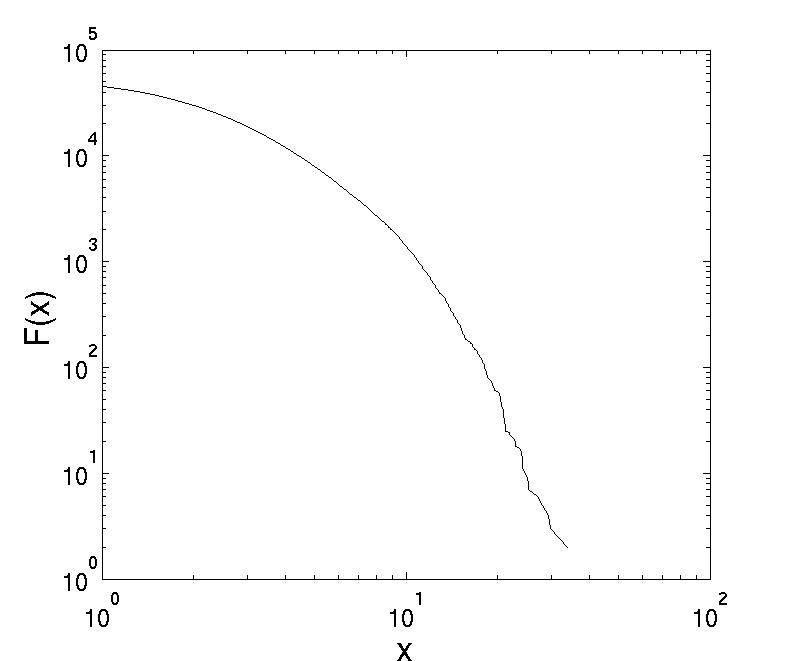}}
\caption{Log-log plots of $F(x)$ for 3D Brownian trees}
\label{fig:DLARectanglePlots}
\end{figure}

\subsection{Brownian Trees.} Unlike in the previous cases, the function $F\paren{x}$ for Brownian trees does not come close to following a power law, so we conclude that they are not \PH statistically self-similar. The plots for one- and two-dimensional homology for 3D Brownian trees are shown in figures~\ref{fig:BT3D1Plot} and~\ref{fig:BT3D2Plot} (a plot for 2D Brownian trees is similar). The data were gathered from three Brownian trees of order 30,000. Note that the plot is significantly bowed out, indicating that the function on a log-log plot is more concave than linear. This result is expected, as Mandelbrot et. al. have shown with high statistical certainty that Brownian trees are not self-similar~\cite{crosscut}~\cite{deviations}.

\begin{figure}[htp]
\centering
\includegraphics[width=\textwidth]{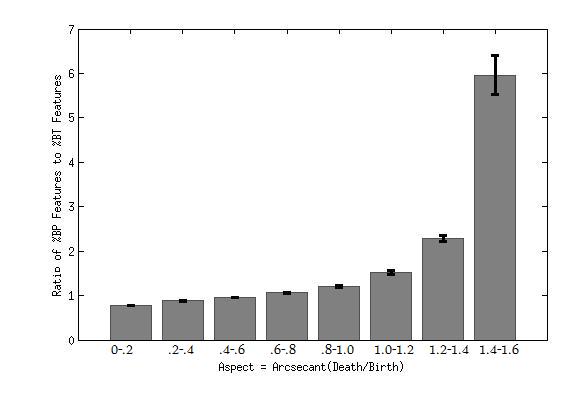}
\caption{
Ratio of branched polymer features to Brownian tree features for various aspects}
\label{fig:Chart}
\end{figure}

\subsection{Shape Measures Involving the Aspect $y$.} The function $F\paren{x}$ is not the only way to glean useful information from the Persistent Homology of polymers.

Consider the bar chart in Figure~\ref{fig:Chart}. For each interval $I$ of possible aspects listed on the horizontal axis($0-.2$, $.2-4$, etc.), the height of the bar is the ratio $P_{BP}(I)/P_{BT}(I)$, where $P_{BP}(I)$ is the proportion of two-dimensional branched polymer features ($i=1$) whose aspect $y$ is in the interval $I$ and $P_{BT}(I)$ is the proportion of two-dimensional Brownian tree features ($i=1$) whose aspect $y$ is in the interval $I$. Recall that the aspect $y$ of a persistent homology class is equal to $\,\arcsec\paren{{d}/{b}},$ where $d$ and $b$ are the death point and birth point of the class, respectively. This effectively measures the angular opening at the edge of a gulf in the structure.

The proportion of features with low aspect measurements appearing in the branched polymers and the Brownian trees is similar. However, features with aspect greater than $1.2$ are more common in branched polymers and those with aspect greater than $1.4$ are much more abundant. This makes sense when one examines the pictures of the branched polymer and the Brownian tree shown in Figure~\ref{fig:polymers}. As described in \S\ref{Alexander Duality}, features with a large aspect value are gulfs  with the property that a large ball can fit inside them but only a very small ball can enter them. Intuitively, this difference seems to capture the fact that most features of the Brownian trees have a similar shape. These gaps between the branches of the Brownian tree can trap balls only slightly larger than those that can enter them. On the other hand, branched polymers have more varied gulfs.

\section{Conclusions.}  Geometric measures of shape suitable for probabilistically generated fractal structures like branched polymers, Brownian trees, and self-avoiding walks are not so plentiful.  These structures lack differentiability, so curvature measures aren't defined.  They are contractible as topological spaces, so the usual the usual measures of topology (homology, homotopy, etc.) are trivial.

The topology of  the $\epsilon$-neighborhoods $S_\epsilon$ of $S$  contains more useful geometric information about $S$.  Via persistent homology, the topology of these neighborhoods gives rise to a \PH density function $f(x,y)$, for each degree $i$, that depends on two variables, the size $x$ and the aspect $y$.  The \PH density function is a subtle measure of shape.

Our theoretical observations are:
\begin{itemize}
\item This measure of shape can be interpreted as how the set $S$ can ``trap'' $\epsilon$-thickened objects of dimension $m-i-1$, e.g. $\epsilon$-balls if $i=m-1$.
\item The $x$ dependency of $f$ gives rise to a notion of fractal dimension, the {\em \PH dimension}.
\item For some examples like the generalized Sierpinski triangle, the \PH  dimension is independent of $i$, and agrees with the Hausdorff dimension.
\end{itemize}

Our main results are computational, concerning branched polymers for $m=2$ or $3$, Brownian trees for $m=2$ or $3$, and self-avoiding walk for $m=2$.  We show that:
\begin{itemize}
\item Significant quantitative information about $f(x,y)$ can be obtained by calculating scatter plots of \PH points for individual examples.
\item For branched polymers, $m=3$, the \PH dimension is about $2$ for both $i=1$ and $i=2$, which agrees with the known estimate for Hausdorff dimension.
\item For self-avoiding walk, $m=2$ and $i=1$, the \PH dimension is approximately $4/3$, which agrees with its conjectured equality with the Hausdorff dimension of SLE~$8/3$.
\item For Brownian trees, the function $F(x)$ is strictly concave on a log-log plot, so Brownian trees are not \PH statistically self-similar.
\item The $y$ dependence of $f(x,y)$ is significantly different for branched polymers and Brownian trees, which quantifies visible differences between the shapes of the two structures.
\end{itemize}

This paper is intended as a feasibility demonstration for a method, rather than as a conclusive compilation of results.  We hope that the shape measure given by the \PH density function will be useful in other situations with complicated geometric structures that are probabilistically determined.  Examples we would like to see include spatial distribution of turbulence,  defect structures in materials, foams, and many others.

\vspace{1em}\noindent{\bf Acknowledgments.}
We would like to thank Tom Spencer and Jeremy Mason for useful conversations, and Benjamin Mann for his enthusiastic support.  The Institute for Advanced Study and DARPA provided support for this project.

\newcommand{\quash}[1]{}
\quash{

\section{Stuff to Add}
(Persistent Homology Dimension)
There is no inequality between  and Hausdorff dimension.
PHD measures complexity.
homologically self-similar -- comes with a similarity ratio
Homologically SSS all ratios

No relation between PH self-similar and self-similar.

independent of embedding and embedding dimension

reduced homology

connected case i=0

all deaths before diameter

Add edelsbrunner to the Ghrist citation

mouth of the gulf

Two length scales for our examples

}
\clearpage

\bibliographystyle{plain}
\bibliography{HomBib}

\end{document}